\documentclass[11pt]{amsart}
\usepackage{amsmath}
\usepackage{amssymb}
\usepackage{graphicx}
\usepackage[utf8]{inputenc}

\newtheorem{theorem}{Theorem}[section]

\newtheorem{lemma}[theorem]{Lemma}

\theoremstyle{definition}
\newtheorem{definition}[theorem]{Definition}

\newtheorem{remark}[theorem]{Remark}
\newtheorem{acknowledgement}[theorem]{Acknowledgement}
\numberwithin{equation}{section} \email{{\tt jidiaz@ucm.es}}
\input{tcilatex}

\begin{document}
\title[New applications of monotonicity methods]{New applications of
monotonicity methods to a class of non-monotone parabolic quasilinear
sub-homogeneous problems}
\author[J.I. D\'{\i}az]{Jes\'{u}s Ildefonso D\'{\i}az}
\address[J.I. D\'{\i}az]{Instituto de Matem\'{a}tica Interdisciplinar,
Universidad Complutense de Madrid, Plaza de Ciencias 3, 28040 Madrid, Spain}

\begin{abstract}
The main goal of this survey is to show how some monotonicity methods
related with the subdifferential of suitable convex functions lead to new
and unexpected results showing the continuous and monotone dependence of
solutions with the respect to the data (and coefficients) of the problem. In
this way, this paper offers `a common roof' to several methods and results
concerning monotone and non-monotone frameworks. Besides to present here
some new results, this paper offers also a peculiar review to some topics
which attracted the attention of many specialists in elliptic and parabolic
nonlinear partial differential equations in the last years under the
important influence of Ha\"{\i}m Brezis. To be more precise, the model
problem under consideration concerns to positive solutions of a class of
doubly nonlinear diffusion parabolic equations with some sub-homogeneous
non-monotone forcing terms.
\end{abstract}

\maketitle

\section{Introduction}

This survey offers \emph{a common roof} to several methods and results
concerning the continuous dependence of solutions with respect to the data
in monotone and non-monotone frameworks. So, besides to present here some
new results, this paper offers also a peculiar survey to some topics which
attracted the attention of many specialists in elliptic and parabolic
nonlinear partial differential equations in the last years. We will show how
some monotonicity methods (as in Brezis \cite{Brezis Monotonic} and Lions
\cite{Lions} ), related with the subdifferential of suitable convex
functions, lead to new results concerning the monotone and continuous
dependence of solutions on an unexpected framework for the problem under
consideration. Our main goal here is not exactly the existence of solutions
but the continuous and monotone dependence of solutions with respect to the
data (and coefficients) of the problem in $L^{2}$ when the expected space
for it is reduced to $L^{1}$. Most of the result of this paper will deal
with positive solutions of the following class of doubly nonlinear diffusion
parabolic equations (in divergence form) with a sub-homogeneous non-monotone
forcing term
\begin{equation*}
(P)\qquad \left\{
\begin{array}{rclr}
\partial _{t}(u^{2q-1})-\Delta _{p}u & = & f(x,u)+h(t,x)u^{q-1} & \text{ in }%
Q_{T}\overset{}{:=}(0,T)\times \Omega , \\[0.1cm]
u & = & 0\  & \text{ on }\Sigma :=(0,T)\times \partial \Omega , \\[0.1cm]
u(0,.) & = & u_{0}(.) & \text{ on }\Omega ,%
\end{array}%
\right.
\end{equation*}%
where $\Omega $ is a smooth bounded domain in $\mathbb{R}^{N}$, $N\geq 1$, $%
T>0$ and with $\Delta _{p}u$ the usual $p-$Laplacian operator, $\Delta _{p}u=%
\mathrm{div}(|\nabla u|^{p-2}\nabla u)$ for $1<p<\infty $. We emphasize that
probably the interest of our results is not for the applications to the
above doubly nonlinear equations but by its method of proof. Moreover, they
are new even for the case of a linear diffusion as $(P)$ with $p=2$. We
assume in $(P)$ a possible nonlinear inertia term (i.e. in the time
derivative), for some
\begin{equation}
q\in (1,p]  \label{Hypo q}
\end{equation}%
and a \textit{sub-homogeneous} forcing term $f(x,u)+h(t,x)u^{q-1}$, where
\begin{equation}
h\in L^{1}(0,T:L^{2}(\Omega )),  \label{Hypo L^1 h(x,t)}
\end{equation}%
and with the non-homogeneous perturbation term $f(x,u)$ satisfying the
following structural assumptions:

\begin{itemize}
\item[(f1)] \label{f1}$f(x,u)$ is a continuous function on $u\in (0,+\infty
) $, for a.e. $x\in \Omega $ and $x\rightarrow f(x,u)$ belongs to $%
L^{2}(\Omega ),$ for any $u\in (0,+\infty )$,

\item[(f2)] \label{f2}$f(x,u)=f_{1}(x,u)+f_{2}(x,u)$ with $\dfrac{f_{1}(x,u)%
}{u^{q-1}}$ non increasing and $\dfrac{f_{2}(x,u)}{u^{q-1}}$ is globally
Lipschitz continuous in $u\in (0,+\infty )$, of Lipschitz constant $K\geq 0,$
for a.e. $x\in \Omega $,

\item[(f3)] $\lim_{r\downarrow 0}\dfrac{f_{1}(x,r)}{r^{q-1}}=a_{0}(x)$ with $%
a_{0}\in L^{2}(\Omega )$.
\end{itemize}

\noindent Additionally, in some cases, we shall need also the condition

\begin{itemize}
\item[(f4)] for any $z>0$ there exists $v_{z}\in L^{\infty }(\Omega )$ such
that $z=\dfrac{f_{1}(x,v_{z}(x))}{\left\vert v_{z}(x)\right\vert ^{q-1}}%
-a_{0}(x)$ a.e. $x\in \Omega $.
\end{itemize}

\noindent Notice that, since we shall not pay attention to the existence of
solutions but to the continuous dependence with respect to the data, \textit{%
no sign condition }is assumed on $h(t,x)$ although we are interested in
positive solutions of $(P)$. Notice also that, as in \cite{Di Thelin},
condition (f2) can be simply formulated as
\begin{equation*}
f(x,u)-f(x,\widehat{u})\geq -K\big (u^{q-1}-\widehat{u}^{q-1}\big )\text{
for any }u>\widehat{u}\geq 0\text{ and a.e. }x\in \Omega .
\end{equation*}%
Condition (f4), of technical nature, will be required only when $f_{1}(x,r)$
is $x-$dependent and express some kind of surjectivity condition of the
application $u\mapsto \frac{f_{1}(x,u)}{u^{q-1}}$, over $(0,+\infty ).$ We
also point out that assumptions (f1) and (f4), for some $q\in (1,p]$, are
compatible with other assumptions, near $r=0$ and near $r=+\infty $, which
arise in the literature and that allows to consider some singular problems.
For instance, in \cite{Diaz-Saa} it was proved that the necessary and
sufficient condition for the existence of a positive solution for the
stationary problem associated to $(P),$ when $h(t,x)=K=0$ is that
\begin{equation*}
\lambda _{1}(-\Delta _{p}v-a_{0}v^{p-1})<0
\end{equation*}%
and%
\begin{equation*}
\lambda _{1}(-\Delta _{p}v-a_{\infty }v^{p-1})>0,\text{ }a_{\infty
}(x)=\lim_{r\uparrow +\infty }\frac{f_{1}(x,r)}{r^{p-1}}.
\end{equation*}%
There are many variants in the literature: for instance, in \cite%
{Ghergu-Radulescu libro poblacion} (see page 275) it is assumed (for $p=2$)
that $\lim_{r\downarrow 0}\frac{f_{1}(x,r)}{r^{p-1}}=+\infty $ and that $%
\lim_{r\uparrow +\infty }\frac{f_{1}(x,r)}{r^{p-1}}=0$.

On the initial condition we will assume that
\begin{equation}
u_{0}\in L^{2q}(\Omega )\cap W_{0}^{1,p}(\Omega ),\text{ }u_{0}>0\text{ on }%
\Omega .\text{ }  \label{Hypo u_0}
\end{equation}%
but some more general conditions are also possible (see Remark \ref{Rem Hipo
dato inicial}).

Very often the nonlinear diffusion equation is equivalently written, in
terms of $W=u^{2q-1}$ with
\begin{equation*}
m=\frac{1}{2q-1}\in \lbrack \frac{1}{2p-1},1)
\end{equation*}%
as
\begin{equation*}
(P_{m,p,q})\qquad \left\{
\begin{array}{rclr}
\partial _{t}W-\Delta _{p}W^{m} & = & f(x,W^{m})+h(t,x)(W^{m})^{\frac{1-2m}{%
2m}} & \text{ in }Q_{T}, \\
W^{m} & = & 0\  & \text{ on }\Sigma , \\
W(0,.) & = & u_{0}^{2q-1} & \text{ on }\Omega .%
\end{array}%
\right.
\end{equation*}%
Since $(p-1)m=\frac{p-1}{2q-1}\in \lbrack \frac{p-1}{2q-1},p-1)$, the
diffusion operator in problem $(P)$, i.e. $(P_{m,p,q})$, offers three
different classes of diffusions, in the terminology of \cite{Di Herrero},
\cite{Kalashnikov}, \cite{Di}, \cite{Tsutsumi}, \cite{Gilding-Kersner}, \cite%
{Dia Extrac}, \cite{Vazquez libro}:

i) \emph{fast diffusion} (which corresponds to $(p-1)m<1$, $i.e.$ $q\in
(\max (\frac{p}{2},1),p]$),

ii) \emph{slow diffusion} (which corresponds to $(p-1)m>1$, $i.e.$ $p>2$ and
$q\in (1,\frac{p}{2})$),

and

iii) the case $(p-1)m=1$ ($i.e.$ $q=\frac{p}{2})$, which was considered, for
instance, in \cite{Del Pino-Daulbaut} in connection with \textit{optimal
logarithmic Sobolev inequalities}: see also \cite{Saa}.

Since the perturbation in the right hand side can be written as $(W^{m})^{%
\frac{1-2m}{2m}}=W^{r}$ with $r:=\frac{1-2m}{2},$ if we assume, for
instance, $p=2$ then $m\in \lbrack \frac{1}{3},1)$ and, in particular $%
0<r<m<1$: a case considered for $h=1$ and $f=0$ by several authors as, e.g.
\cite{Pablo-Vazquez}, and \cite{Hui}: see also \cite{Guo}.

In the limit case $q=1$ ($i.e.$ $m=1$), the problem \textit{formally}
includes a Heaviside function (a model similar to the one which appears in
some climate models with the $p$-Laplace operator) since, roughly speaking,
we can approximate the problem by other ones corresponding to a sequence of
exponents $q_{n}$ $\searrow 1$ as $n\rightarrow +\infty $ and thus it seems
possible to extend the conclusions to the multivalued problem%
\begin{equation*}
(P_{H})\qquad \left\{
\begin{array}{rclr}
\partial _{t}W-\Delta _{p}W & \in & f(x,W)+h(t,x)H(W) & \text{ in }Q_{T}, \\%
[0.1cm]
W & = & 0\  & \text{ on }\Sigma , \\[0.1cm]
W(0,.) & = & W_{0} & \text{ on }\Omega ,%
\end{array}%
\right.
\end{equation*}%
with $H(r),$ the Heaviside, multivalued-function, $H(r)=\{0\}$ if $r<0,$ $%
H(r)=\{1\}$ if $r>0$ and $H(0)=[0,1]$. Problems similar to $(P_{H})$ appear
in many contexts, and, in particular, in climate Energy Balance Models (see,
e.g., \cite{Di-Tello}, \cite{Bensid-Di}, and their references). \ For some
comparison results concerning solutions of $(P_{m,p,q})$ corresponding to
two different values of $m$ see \cite{Benilan-Diaz}. The continuous
dependence on $m$ (even in a more general framework than the one here
considered) was studied in \cite{Beni-Cr depende} and \cite{Beni-Cr-Saks}.

It is well known (see, $e.g.$, the exposition made in \cite{Brezis Monotonic}%
, \cite{Ben Madrid}, \cite{Di Thelin}, \cite{Dia Extrac}) that the theory of
maximal monotone operators on Hilbert spaces [or, more in general, the
theory of m-accretive operators in Banach spaces: see, $e.g.$, \cite{Barbu},
\cite{BenCran Pazy} and the surveys \cite{Evans} and \cite{Benilan-Wight}]
can be applied to the above class of problems in the absence of the forcing
term or when it is assumed to be globally Lipschitz continuous on the
corresponding functional space. But it seems that the applicability of the
abstract theory of such type of operators is not well known in the
literature when the forcing term is merely \textit{sublinear} (if $p=2$) or,
more generally, \textit{sub-homogeneous }($q\leq p$ if $p\neq 2$). For some
pioneering results we send the reader to \cite{Kranoselki}, \cite{Fujita
-Watananbe}, \cite{Keller-Cohen}, \cite{Keller Krisne}, \cite{Amman}, \cite%
{Berestycki}, \cite{Levine-Sacks}, \cite{Levine SIAM} and the book \cite%
{Samarski-Galaktionov libro}.

As said before, the main goal of this paper is to show how the above
mentioned \emph{monotonicity methods} can be suitably applied also to this
class of non-monotone problems, leading to a general framework (specially
concerning the $x$-dependence of coefficients) in which it is possible to
show the \textit{continuous and monotone dependence with respect to the data}
(the initial datum and the \emph{potential type} coefficient $h(t,x)$) even
if there are non-monotone terms in the right hand side.

As a matter of fact, in contrast with the previous literature, we will show
that it is possible to give a sense to the solvability of the equation even
for time dependent coefficients $h(x,t)$ satisfying merely (\ref{Hypo L^1
h(x,t)}) (see some comments on the difficulties arising when using a more
classical variational approach in \cite{Boccardo-Orsina}, \cite{Artola},
\cite{Palmieri}) and, what it is more important, \textit{without prescribing
any sign} on $h(x,t)$, which corresponds to the so-called \textit{indefinite
perturbed problems} arising, for instance, in population dynamics: see \cite%
{Namba}, \cite{Bandle-Pozio-Tesei}, \cite{Badii-Diaz-Tesei} and \cite{Arias
Cuesta Gossez}, among many other possible references.

As we will see, it is useful to start our program by considering the\textit{%
\ }sub-homogeneous \textit{simpler }problem corresponding to $f(x,u)\equiv 0,
$ $i.e.$ the problem
\begin{equation*}
(P_{q})\qquad \left\{
\begin{array}{rclr}
\partial _{t}(u^{2q-1})-\Delta _{p}u & = & h(t,x)u^{q-1} & \text{ in }Q_{T},
\\
u & = & 0\  & \text{ on }\Sigma , \\
u(0,.) & = & u_{0}(.) & \text{ on }\Omega .%
\end{array}%
\right.
\end{equation*}%
The existence and uniqueness of a $L^{1}-$mild positive solution when $%
h(t,x)\leq 0$ is a consequence of the well-known m-T-accretivity results of
the associated operator (see, e.g. \cite{Benilan Franco-Japan}, \cite{Di
Thelin} and\ \cite{Vazquez libro}). Nevertheless, since the right hand side
is non-Lipschiz continuous, problem $(P_{q})$ (and also problem $(P)$) may
have more than one solution (in particular when $h(t,x)$ is changing sign
and negative near $\Sigma $ and we assume $p>2$ and $q\in (1,\frac{p}{2})$)$.
$ Nevertheless we can introduce a method to select only one $L^{1}-$mild
positive solution by means of some monotonicity arguments. Indeed, we will
select the $L^{1}-$mild positive solution $u$ of $(P_{q})$ such that $w(%
\frac{q}{2q-1}t)=u(t)^{q}$ coincides with the unique $L^{2}-$mild positive
solution of the problem%
\begin{equation}
\left\{
\begin{array}{lc}
\dfrac{dw}{dt}+\partial J_{0,q}(w)\ni h(t) & \text{in }L^{2}(\Omega ), \\%
[0.2cm]
w(0)=w_{0}, &
\end{array}%
\right.   \label{Problem parabolic subdifferential J-0q}
\end{equation}%
where $J_{0,q}$ is the functional in $L^{2}(\Omega )$ given by%
\begin{equation*}
J_{0,q}(w)=\left\{
\begin{array}{lr}
\displaystyle\dfrac{q}{p}\int_{\Omega }|\nabla w^{\frac{1}{q}}|^{p}\mathrm{d}%
x & \text{if }w\in D(J_{0,q}), \\
+\infty  & \text{otherwise,}%
\end{array}%
\right.
\end{equation*}%
with
\begin{equation*}
D(J_{0,q}):=\{w\in L^{2}(\Omega )\text{ such that }w\geq 0\text{ and }w^{%
\frac{1}{q}}\in W_{0}^{1,p}(\Omega )\}.
\end{equation*}%
Developing an idea of D\'{\i}az and Sa\'{a} \cite{Diaz-Saa} (for $p\neq 2$)
we will see that $J_{0,q}$ is a convex, lower semicontinuous functional and
thus its subdifferential $\partial J_{0,q}(w)$ is well defined and the
uniqueness of a $L^{2}-$mild positive solution $w$ of (\ref{Problem
parabolic subdifferential J-0q}) is well-known. In that case we say that $%
u(t)$ is the \textit{selected }$L^{1}-$\textit{mild positive solution} of $%
(P_{q})$ (and so it is unique). Of course that if, under some additional
assumptions, it can be shown the uniqueness of a positive weak solution of
the equation then necessarily it must coincides with the selected $L^{1}-$%
mild positive solution (see, e.g., \cite{Pablo-Vazquez}, \cite{Hui}, \cite%
{Guo}, \cite{Cazenave-Dik-Escobedo}, \cite{Dickstein} and \cite{Charro-Peral}%
, among others).

In Section 2 of this paper we will study the subdifferential $\partial
J_{0,q}(w)$. We we will prove that, given $\mu >0$ and $h\in L^{2}(\Omega )$%
, the resolvent equation%
\begin{equation}
w+\mu \partial J_{0,q}(w)\ni h  \label{Resolvent equation 2}
\end{equation}%
is connected, through the relation $w=u^{q}$, with the auxiliary variational
problem
\begin{equation*}
\underset{v\in K}{\min }J_{h,q}(v)
\end{equation*}%
where%
\begin{equation*}
K:=\big \{v\in W_{0}^{1,p}(\Omega )\cap L^{2q}(\Omega ),v\geq 0\text{ on }%
\Omega \big \}
\end{equation*}%
and%
\begin{equation*}
J_{h,q}(v):=\frac{\mu }{p}\int_{\Omega }|\nabla v|^{p}\mathrm{d}x+\frac{1}{2q%
}\int_{\Omega }|v|^{2q}\mathrm{d}x-\frac{1}{q}\int_{\Omega }h(x)|v|^{q}%
\mathrm{d}x.
\end{equation*}

Since the problem is sub-homogeneous ($q\in (1,p]$) the different terms of $%
J_{h,q}(v)$ satisfy good growth conditions and the existence and uniqueness
of a minimum $v_{h,q}\in K$ can be obtained by standard direct methods of
the Calculus of Variations (see, $e.g.$, Lemma 5 of \cite{Benguria Brezis
Lieb} for the case $p=2$ and \cite{Takac-Tello-Ulm} for $p>1$ and $q\in
\lbrack 1,p]$). Once again, the Euler-Lagrange equation
\begin{equation}
-\mu \Delta _{p}v+v^{2q-1}=h(x)v^{q-1}\text{ in }\Omega ,
\label{Euler-Lagrange}
\end{equation}%
may have other weak solutions $\widehat{v}\in W_{0}^{1,p}(\Omega )$
different to the minimum $v$ of $J_{h,q}$ (specially if the sign of $h(x)$
is not prescribed, $h(x)$ is negative near $\Sigma $ and we assume $p>2$ and
$q\in (1,\frac{p}{2})$)) but the relation $w=u^{q}$ allows to select only $v$
when we assume that $w$ is the solution of (\ref{Resolvent equation 2}).

As we shall show in Section 3, the definition of a unique $u(t)$ \textit{%
selected }$L^{1}-$\textit{mild positive solution} of $(P)$ can be also
obtained for the general case of $f\neq 0$ as indicated before by following
a similar process to the indicated above. The main result of this paper is
the following:

\begin{theorem}
\label{Theore main} Let $q\in (1,p]$ and $h\in L^{1}(0,T:L^{2}(\Omega ))$.
Let $u_{0},f$ satisfying (\ref{Hypo u_0}) and (f1)-(f3). Assume that $%
f_{1}(x,u)=f_{1}(u)$ independent of $x$ or $f_{1}(x,u)$ satisfying also
(f4). Then for any $T>0$, there exists a unique selected positive $L^{1}-$%
mild solution $u$ to problem $(\mathrm{P})$ and $u^{q}\in
C([0,T]:L^{2}(\Omega ))$. In addition, if $h\in L^{\infty }(0,T:L^{\infty
}(\Omega ))$ and $u_{0}\in L^{\infty }(\Omega )$ then $u\in L^{\infty
}(0,T:L^{\infty }(\Omega ))$. Moreover, if $v_{0}$ and $g$ satisfy the same
conditions than $u_{0}$ and $h$, and if $v$ is the respective selected
positive $L^{1}-$mild solution of problem $(\mathrm{P}),$ then, for any $%
t\in \lbrack 0,T]$ we have the monotone continuous dependence estimate%
\begin{equation}
\begin{array}{ll}
\Vert (u^{q}(t)-v^{q}(t))^{+}\Vert _{L^{2}(\Omega )} & \leq e^{Kt}\Vert
(u_{0}^{q}-v_{0}^{q})^{+}\Vert _{L^{2}(\Omega )} \\[0.1cm]
& \displaystyle+\int_{0}^{t}e^{K(t-s)}\left\Vert [h(s)-g(s)]_{+}\right\Vert
_{L^{2}(\Omega )}ds,%
\end{array}
\label{Estimate  dependence parabolic}
\end{equation}

where $K\geq 0$ is the constant indicated in (f2).
\end{theorem}

\bigskip

Notice that, in particular, for the case of a slow diffusion, $p>2$ and $%
q\in (1,\frac{p}{2})$, the above conclusions hold for `flat solutions' ($%
i.e. $ positive solutions such that $u=\frac{\partial u}{\partial n}=0$ on $%
\Sigma $). Notice that even for the special case $h=g$ estimate (\ref%
{Estimate dependence parabolic}) is new for the doubly nonlinear problem $(%
\mathrm{P})$: indeed, as indicated before the accretivity results of the
doubly nonlinear diffusion operator leads only to $L^{1}(\Omega )-$monotone
continuous dependence estimates (if $p=2$ such estimates also hold on $%
H^{-1}(\Omega )$ \cite{Brezis Monotonic}$)$, but not in $L^{2}(\Omega )$
(see, e.g., B\'{e}nilan \cite{Benilan not L`2}) as it is expressed in (\ref%
{Estimate dependence parabolic}).

We point out that, obviously, the function $u_{\infty }(x)\equiv 0$ in $%
\Omega $ is a trivial solution of the stationary problem associated to $(P)$%
. Here we are interested on positive solutions of problem $(P)$. We will
prove (see Theorem \ref{Theorem extinction copy(1)}) that, in fact, if $q\in
(1,\frac{p}{2})$ and $p>2,f(x,u)\equiv 0,~h\in L^{1}(0,T:L^{2}(\Omega )),$ $%
h\geq 0$ and $u_{0}\gvertneqq 0$ then there is no extinction in finite time,
so that $\Vert u^{2q-1}(t)\Vert _{L^{2}(\Omega )}>0$ for any $t>0.$ The
situation is different if $q\in (\frac{p}{2},p]$ since, at least for $%
f(x,u)\equiv 0$ and $h\leq 0$, there is a finite extinction time $T_{e}>0$,
such that $w(t)\equiv 0,$ in $\Omega $, for any $t\geq T_{e}$. In that case,
we understand that the $L^{1}-$mild solution $u(t)$ of $(P)$ also
extinguishes in $\Omega $ after $T_{e}.$

In the Section 3 we will study of the auxiliary simplified problem $(P_{q})$
through the study of the subdifferential operator $\partial J_{0,q}(v)$ in $%
L^{2}(\Omega )$. This will allow to get the proof of Theorem \ref{Theore
main} by application of some abstract results on monotone operators on
Hilbert spaces. Many other variants, commented in form of a series of
Remarks, opening the application of this view point to many other different
formulations, will be presented. This is the case, for instance when the $p$%
-Laplacian is replaced by an homogeneous diffusion operator of the form $%
\mathrm{div}(a(x,\nabla u))$ with the homogeneity condition
\begin{equation*}
A(x,t\mathbf{\xi )=}\left\vert t\right\vert ^{p}A(x,\mathbf{\xi )}\text{ for
all }t\in \mathbb{R}\text{ and all }(x,\mathbf{\xi )\in }\Omega \times
\mathbb{R}^{N},
\end{equation*}%
where $a(x,\mathbf{\xi })=\frac{1}{p}\partial _{\mathbf{\xi }}A(x,\mathbf{%
\xi )}$.

\section{On the subdifferential of $J_{0,q}$}

The proof of the main results will be obtained through the study of the
Cauchy problem
\begin{equation*}
\left\{
\begin{array}{lc}
\dfrac{dw}{dt}+\partial J_{0,q}(w)\ni h(t) & \text{in }L^{2}(\Omega ) \\%
[0.2cm]
w(0)=w_{0}, &
\end{array}%
\right.
\end{equation*}%
with $J_{0,q}$ the functional presented in the Introduction. The convexity
of $J_{0,q}$ will play a crucial role in the rest of the paper.

\bigskip

\begin{lemma}
\label{Convexity} Given $q\in (1,p]$, the functional $J_{0,q}$ is convex,
lower semicontinuous and proper on $L^{2}(\Omega )$.
\end{lemma}

\noindent \textit{Proof}. The proof for the case $q=p$ was given in Lemma 1
of \cite{Diaz-Saa}, and the proof for the case $q\in (1,p)$ was obtained in
\cite{Takac-Tello-Ulm} (see Lemma 4 and Example 5.2). A different proof of
this last case can be obtained from Proposition 2.6 of \cite{Brasco-Franzina}%
. To prove that $J_{0,q}$ is lower semicontinuous in $L^{2}(\Omega )$ it
suffices to prove that if we have a sequence $\rho _{n}\rightarrow \rho $ in
$L^{2}(\Omega )$ such that $J_{0,q}(\rho _{n})\leq \lambda $ then $%
J_{0,q}(\rho )\leq \lambda $. But since $\rho _{n}^{1/q}$ is bounded in $%
W_{0}^{1,p}(\Omega )$ there exists a subsequence, still labeled as $\rho
_{n}^{1/q}$, such that $\rho _{n}^{1/q}$ converges weakly in $%
W_{0}^{1,p}(\Omega )$, so that $\nabla \rho _{n}^{1/q}$ converges weakly in $%
L^{p}(\Omega )^{N}$ and since the norm is lower semicontinuous we obtain
that $\liminf_{n}J_{0,q}(\rho _{n})\geq J_{0,q}(\rho )$, and hence $%
J_{0,q}(\rho )\leq \lambda ._{\quad \blacksquare }$

\bigskip

\begin{remark}
As indicated in \cite{Diaz-Saa}, the main results of \cite{Diaz-Saa} were
presented in September 1985 in \cite{Diaz-Saa CEDYA}. Its Lemma 1 extends
and develops to the case $p\neq 2$ Remark 2 of Brezis and Oswald \cite{Bre
Oswald} which was inspired in the paper Benguria, Brezis and Lieb \cite%
{Benguria Brezis Lieb} where some previous results of Rafael Benguria's
Ph.D. thesis \cite{Beguria tesis} were presented together with some newer
results. So, in contrast to what is indicated in \cite{Brasco-Franzina}, the
consideration of the case $p\neq 2$ was not carried for the first time in
\cite{Belloni-Kawhol} but in \cite{Diaz-Saa CEDYA}, \cite{Diaz-Saa}
seventeen years before. The extension to the case of $\mathbb{R}^{N}$ was
carried out in \cite{Chaib} (for an extension to weaker solutions see \cite%
{Dat}).
\end{remark}

\bigskip

\begin{remark}
It seems, that the connection between Lemma 1 of \cite{Diaz-Saa} (called by
some authors \emph{D\'{\i}az -Sa\'{a} inequality} when $q=p,$ \cite{Chaib},
\cite{Takac-Giacomoni}) and the generalization of the 1910 Picone inequality
\cite{Picone} (concerning originally with ordinary differential equations
and much more later extended to partial differential equations in \cite%
{Allegreto-Huang}; see, also the survey \cite{Dosly}) was pointed out for
the first time in Chaib \cite{Chaib}. As a matter of fact, it was proved in
Section 3.2 of \cite{Brasco-Franzina} that the convexity of $J_{0,q}$ (for
any $q\in (1,p]$) is equivalent to the generalized Picone inequality
\begin{equation*}
\frac{1}{p}\left\vert \nabla u\right\vert ^{p-2}\langle \nabla u,\nabla
\left( \frac{z^{q}}{u^{q-1}}\right) \rangle \leq \frac{q}{p}\left\vert
\nabla z\right\vert ^{p}+\frac{p-q}{p}\left\vert \nabla u\right\vert ^{p}%
\text{ a.e. on }\Omega
\end{equation*}%
if $u,z\in W_{loc}^{1,p}(\Omega ),~u>0,~z\geq 0\text{ on }\Omega $.
\end{remark}

\bigskip

We recall that given a convex, l.s.c. function $\phi :H\rightarrow (-\infty
,+\infty ]$, $\phi $ proper, over a Hilbert space $H$, a pair $(w,z)\in
H\times H$ is such that $z\in \partial \phi (w)$ if $\forall \xi \in H,$ $%
\phi (\xi )\geq \phi (w)+(z,\xi -w)$. We say that $w\in D(\phi ):=\{v\in H$
such that $\phi (v)<+\infty \}$ is such that $w\in D(\partial \phi )$ if the
set of $z\in \partial \phi (w)$ is not empty. We have

\begin{equation*}
D(\partial J_{0,q})\subset D(J_{0,q})\subset \overline{D(J_{0,q})}^{L^{2}}=%
\overline{D(\partial J_{0,q})}^{L^{2}}
\end{equation*}%
(see Proposition 2.11 of Brezis \cite{Brez OMXM}). The following result
proves that the operator $\partial J_{0,q}$ satisfies an additional property
to the mere monotonicity: it is a $T-$monotone operator in $L^{2}(\Omega )$
in the sense of Brezis-Stampacchia (\cite{Brezis-Stampacchia}). This will
explain later the comparison of solutions of problem $(P)$ with respect to
different data $h(t,x)$ for solutions.

\begin{lemma}
\label{Prop Benilan} Let $\tau (s)=s_{+}$. Then for any $w,\widehat{w}\in
L^{2}(\Omega )$%
\begin{equation}
J_{0,q}\left( w-\tau (w-\widehat{w})\right) +J_{0,q}\left( \widehat{w}+\tau
(w-\widehat{w})\right) \leq J_{0,q}(w)+J_{0,q}(\widehat{w}).
\label{tau-contraction}
\end{equation}

\noindent In particular $\partial J_{0,q}$ is a $T-$monotone operator in $%
L^{2}(\Omega )$, i.e. for any $w,\widehat{w}\in D(\partial J_{0,q})$ and $%
z\in \partial J_{0,q}(w),$ $\widehat{z}\in \partial J_{0,q}(\widehat{w}),$%
\begin{equation}
\int_{\Omega }(z-\widehat{z})[w-\widehat{w}]_{+}\mathrm{d}x\geq 0,
\label{T-monotone}
\end{equation}

\noindent and given $h,\widehat{h}\in L^{2}(\Omega )$, if for $\mu >0$, $w,%
\widehat{w}\in L^{2}(\Omega )$ are such that
\begin{equation}
w+\mu \partial J_{0,q}(w)\ni h\text{ and }\widehat{w}+\mu \partial J_{0,q}(%
\widehat{w})\ni \widehat{h},  \label{resolvent equation}
\end{equation}%
then%
\begin{equation}
\left\Vert \lbrack w-\widehat{w}]_{+}\right\Vert _{L^{2}(\Omega )}\leq
\left\Vert \lbrack h-\widehat{h}]_{+}\right\Vert _{L^{2}(\Omega )}.
\label{L^2 contraction}
\end{equation}
\end{lemma}

\noindent \textit{Proof}. Property (\ref{tau-contraction}) is equivalent to
the inequality%
\begin{equation}
J_{0,q}(\min (w,\widehat{w}))+J_{0,q}(\max (w,\widehat{w}))\leq
J_{0,q}(w)+J_{0,q}(\widehat{w}).  \label{max-min contract}
\end{equation}%
Obviously we can assume $w,\widehat{w},\min (w,(\widehat{w}-k)),\max ((w-k),%
\widehat{w})\in D(J_{0,q}):=\{v\geq 0$ and $v^{\frac{1}{q}}\in
W_{0}^{1,p}(\Omega )\cap L^{\frac{2}{q}}(\Omega )\}$ and then, by
Stampacchia's truncation results, we can write%
\begin{equation*}
\int_{\Omega }|\nabla \min (w,\widehat{w})^{\frac{1}{q}}|^{p}\mathrm{d}%
x=\int_{\{w\leq \widehat{w}\}}|\nabla w^{\frac{1}{q}}|^{p}\mathrm{d}%
x+\int_{\{w>\widehat{w}\}}|\nabla \widehat{w}^{\frac{1}{q}}|^{p}\mathrm{d}x
\end{equation*}%
and%
\begin{equation*}
\int_{\Omega }|\nabla \max (w,\widehat{w})^{\frac{1}{q}}|^{p}\mathrm{d}%
x=\int_{\{w>\widehat{w}\}}|\nabla w^{\frac{1}{q}}|^{p}\mathrm{d}%
x+\int_{\{w\leq \widehat{w}\}}|\nabla \widehat{w}^{\frac{1}{q}}|^{p}\mathrm{d%
}x.
\end{equation*}%
Adding both expressions we get inequality (\ref{max-min contract}). To show
that (\ref{max-min contract}) implies that $\partial J_{0,q}$ is a $T-$%
monotone operator in $L^{2}(\Omega )$, i.e. (\ref{T-monotone}) we shall
develop a suggestion made by H. Brezis in Remark 1.10 of \cite{Brezis Probl
Unilateraux}. Since $z\in \partial J_{0,q}(w)$ and $\widehat{z}\in \partial
J_{0,q}(\widehat{w})$ we know that

\begin{equation*}
\begin{array}{cc}
J_{0,q}(v)-J_{0,q}(w)\geq \int_{\Omega }z[v-w]\mathrm{d}x\geq 0 & \text{for
any }v\in L^{2}(\Omega ), \\
J_{0,q}(v)-J_{0,q}(\widehat{w})\geq \int_{\Omega }\widehat{z}[v-\widehat{w}]%
\mathrm{d}x\geq 0 & \text{for any }v\in L^{2}(\Omega ).%
\end{array}%
\end{equation*}%
By taking $v=\min (w,\widehat{w})=w-[w-\widehat{w}]_{+}$ in the first of the
two inequalities, and $v=\max (w,\widehat{w})=\widehat{w}+[w-\widehat{w}%
]_{+} $ in the second one, using that
\begin{equation*}
\min (w,\widehat{w})-w=-[w-\widehat{w}]_{+}\text{ and }\max (w,\widehat{w})-%
\widehat{w}=[w-\widehat{w}]_{+},
\end{equation*}%
by adding the results we get

\begin{equation*}
J_{0,q}(\min (w,\widehat{w}))+J_{0,q}(\max (w,\widehat{w}%
))-J_{0,q}(w)-J_{0,q}(\widehat{w})\leq -\int_{\Omega }(z-\widehat{z})[w-%
\widehat{w}]_{+}\mathrm{d}x,
\end{equation*}%
and thus inequality (\ref{max-min contract}) implies property (\ref%
{T-monotone}). By well-known results (see Section IV.4 of Brezis \cite{Brez
OMXM}) we get conclusion (\ref{L^2 contraction}).$_{\substack{ \blacksquare
\\ }}$

\bigskip

\begin{remark}
\label{Remark correccion k>0}For some convex functionals $J$ a stronger
property than (\ref{tau-contraction}) holds:%
\begin{equation}
J(\min (w,(\widehat{w}-k)))+J(\max ((w-k),\widehat{w}))\leq J(w)+J(\widehat{w%
})  \label{tau-contraction for any k}
\end{equation}%
for any $k>0$. This property is equivalent (\cite{Benilan-Picard}) to
inequality (\ref{tau-contraction}) for any $\tau :\mathbb{R\rightarrow R}$
Lipschitz continuous with $0\leq \tau ^{\prime }\leq 1$ and $\tau (0)=0$ and
for any $k>0.$ This property (\ref{tau-contraction}) implies several
important properties for the realization of the operator $w\rightarrow
\partial J(w)$ over the Banach spaces $L^{s}(\Omega )$, $1\leq s\leq +\infty
$ (see Lemma 3 of \cite{Brez-Strauss} and its generalization in a series of
papers (Th\'{e}or\`{e}me 1.2 and Remark 1.4 of \cite{Benilan-Picard}, \cite%
{Beni tesis}, \cite{Benilan cone}) and (\cite{Benilan-Picard})$.$ Property (%
\ref{tau-contraction}) holds for the class of the, so called, \emph{normal
convex functionals }(see the above mentioned references) but to check it for
the special case of the functional $J_{0,q}$ remains as an open problem
(some partial results can be obtained in this direction: see Remark \ref{Rem
T-accretivity in L^1} ).
\end{remark}

\bigskip

An uneasy task is to identify the operator $\partial J_{0,q}$ involved in
the \textit{resolvent equation} (\ref{resolvent equation}) in terms of the
Euler-Lagrange equation associated to the functional $J_{0,q}$. When trying
to do that directly, using merely the functional $J_{0,q}$, \ we see that,
if we assume that $w>0$ on $\Omega $, given a \textit{direction test
function }$\zeta \in W_{0}^{1,p}(\Omega )\cap L^{2}(\Omega )$ the G\^{a}%
teaux derivative of $J_{0,q}$ in $w$ in the direction $\zeta $ is given
formally by%
\begin{equation}
J_{0,q}^{\prime }(w;\zeta )=-\int_{\Omega }\frac{\Delta _{p}(w^{\frac{1}{q}})%
}{w^{\frac{q-1}{q}}}\zeta \mathrm{d}x.  \label{subdiferential calculus}
\end{equation}%
Thus, at least formally, the convexity of $J_{0,q}$ implies the monotonicity
in $L^{2}(\Omega )$ of its subdifferential and so
\begin{equation}
\int_{\Omega }\left (-\frac{\Delta _{p}(w^{\frac{1}{q}})}{w^{\frac{q-1}{q}}}+%
\frac{\Delta _{p}(\widehat{w}^{\frac{1}{q}})}{\widehat{w}^{\frac{q-1}{q}}}%
\right )(w-\widehat{w})\mathrm{d}x\geq 0.  \label{Monotonia en L^2}
\end{equation}

In \cite{Diaz-Saa} it was shown that expression (\ref{subdiferential
calculus}) is well justified if we assume $w\in {\mathcal{D}}(J_{0,q})$ and $%
w,\Delta _{p}(w^{\frac{1}{q}})\in L^{\infty }(\Omega )$. A different
justification was made in Remark 3.3 of Taka\v{c} \cite{Takac Handbook},
this time under the additional condition that $w>0$ on any compact subset $%
M\subset \Omega ,$
\begin{equation*}
\frac{\Delta _{p}(w^{\frac{1}{q}})}{w^{\frac{q-1}{q}}}\in \mathcal{D}%
^{\prime }(\Omega ),
\end{equation*}%
and $w\in C^{0}(\Omega ).$ Nevertheless, it is possible to get some more
general justifications when instead of analyzing separately $J_{0,q}^{\prime
}(w;\zeta )$ we consider the \textit{resolvent equation} (\ref{resolvent
equation}). The following result is inspired by Lemma 6 of \cite{Benguria
Brezis Lieb} concerning a related problem in which $p=q=2$ and $N=3$.

\begin{lemma}
\label{Lemm subdifferential}Given $q\in (1,p]$, $h\in L^{2}(\Omega )$ and $%
\mu >0$, assume that $w\in D(\partial J_{0,q}),$ $w\geq 0$, satisfies the
resolvent equation (\ref{Resolvent equation 2}). Then function $v:=w^{\frac{1%
}{q}}$ satisfies that $v\in W_{0}^{1,p}(\Omega )\cap L^{2q}(\Omega )$, $%
\Delta _{p}v,h(x)v^{q-1}\in L^{1}(\Omega ),$ $v$ is positive in the sense
that
\begin{equation}
\big |\big \{x\in \Omega :v(x)=0\big \}\big |=0,
\label{Positivity Lemma Subdiffer}
\end{equation}

\noindent and $v$ satisfies the sub-homogeneous equation (\ref%
{Euler-Lagrange}) in the sense of distributions. Moreover,

i) if $1<q<p$ and $0<h_{-}(x)=\max (-h(x),0)\leq C_{h_{-}}$near $\partial
\Omega $
\begin{equation}
v(x)\geq Cd(x,\partial \Omega )^{\frac{p}{p-q}}\text{ a.e. }x\in \Omega ,%
\text{ for some }C>0\text{ dependent of }C_{h_{-}},  \label{Decay h negativa}
\end{equation}

ii) if $h_{-}(x)\equiv 0$ near $\partial \Omega $ and $p>2$ with $q\in (1,%
\frac{p}{2})$ then%
\begin{equation}
v(x)\geq Cd(x,\partial \Omega )^{\frac{p}{p-2q}}\text{ a.e. }x\in \Omega ,%
\text{ for some }C>0\text{ independent on }h,  \label{decay q moderate}
\end{equation}

iii) if $h_{-}(x)\equiv 0$ near $\partial \Omega $ and $q\in \lbrack \frac{p%
}{2},p)$ if $p>2,$ or $q\in (\max (1,\frac{p}{2}),p)$ if $p\leq 2$, then%
\begin{equation}
v(x)\geq Cd(x,\partial \Omega )\text{ a.e. }x\in \Omega ,\text{ for some }C>0%
\text{ independent on }h,
\end{equation}

iv) if $q=p$ then
\begin{equation}
v(x)\geq Cd(x,\partial \Omega )\text{ a.e. }x\in \Omega ,\text{ for some }C>0%
\text{ independent on }h.  \label{estricta positividad}
\end{equation}
\end{lemma}

\noindent \textit{Proof}. Since $D(\partial J_{0,q})\subset D(J_{0,q})$ we
know that $v=w^{\frac{1}{q}}\in W_{0}^{1,p}(\Omega )\cap L^{2q}(\Omega )$.
Moreover $h(x)v^{q-1}\in L^{1}(\Omega )$ since $v\in L^{2q-2}(\Omega )$ and $%
h\in L^{2}(\Omega )$. Therefore the equation (\ref{Euler-Lagrange}) has a
meaning in the sense of distributions. Let $\eta \in $ $\widetilde{D}%
:=W_{0}^{1,p}(\Omega )\cap L^{2q}(\Omega )$ (i.e. without the sign condition
$\eta \geq 0$). Define the functional
\begin{equation*}
J_{h,q}(\eta )=\frac{\mu }{p}\int_{\Omega }|\nabla \eta |^{p}\mathrm{d}x+%
\frac{1}{2q}\int_{\Omega }|\eta |^{2q}\mathrm{d}x-\frac{1}{q}\int_{\Omega
}h(x)|\eta |^{q}\mathrm{d}x.
\end{equation*}%
Therefore, for every $\eta \in $ $\widetilde{D}$%
\begin{equation*}
J_{h,q}(v)\leq J_{h,q}(\eta )
\end{equation*}%
so, $v$ is a minimum of $J_{h,q}$ on $\widetilde{D}$. Now, for $\zeta \in
C_{0}^{\infty }(\Omega )$, using that $d(J_{h,q}(v+\epsilon \zeta
))/d\epsilon $ $=0$ we conclude easily that
\begin{equation*}
\mu \int_{\Omega }|\nabla v|^{p-2}\nabla v\nabla \eta \mathrm{d}%
x+\int_{\Omega }v^{2q-1}\eta \mathrm{d}x=\int_{\Omega }h(x)v^{q}\eta \mathrm{%
d}x,
\end{equation*}%
which proves $v$ satisfies (\ref{Euler-Lagrange}) and $\Delta _{p}v\in
L^{2}(\Omega )$. On the other hand,
\begin{equation*}
-\frac{\Delta _{p}(w^{\frac{1}{q}})}{w^{\frac{q-1}{q}}}=h(x)-w\in
L^{2}(\Omega ),
\end{equation*}%
so, necessarily, $w$ is positive (in the sense of (\ref{Positivity Lemma
Subdiffer})). Moreover, using the decomposition $h(x)=h_{+}(x)-h_{-}(x)$,
with
\begin{equation*}
h_{+}(x)=\max (h(x),0),\text{ }h_{-}(x)=\max (-h(x),0),
\end{equation*}%
we can write (\ref{Euler-Lagrange}) as%
\begin{equation*}
-\mu \Delta _{p}v+v^{2q-1}+h_{-}(x)v^{q-1}=h_{+}(x)v^{q-1}\text{ in }\Omega .
\end{equation*}%
The proof of iii) and v) is consequence of the strong maximum principle (%
\cite{Vazquez Fuerte max}, \cite{Pucci-Serrin}) once that $v\geq 0$ on $%
\Omega ,$ $-\mu \Delta _{p}v+v^{2q-1}+h_{-}(x)v^{q-1}\geq 0$ and since the
zero order terms in the above inequality are super-homogeneous ($2q-1\geq p-1
$ and $h_{-}(x)=0$ near $\partial \Omega $ if $q\in \lbrack \frac{p}{2},p)$).

\noindent To prove i) and ii) notice that in both cases there is a \emph{%
strong absorption with respect to the diffusion} once we write
\begin{equation*}
-\mu \Delta _{p}v+v^{2q-1}+h_{-}(x)v^{q-1}=h_{+}(x)v^{q-1}.
\end{equation*}

\noindent In the case ii), if $h_{-}(x)=0$ on a neighborhood $D_{\delta }$
of $\partial \Omega ,$ with $D_{\delta }=\{x\in \Omega :d(x,\partial \Omega
)\leq \delta \}$, for some $\delta >0$, then $-\mu \Delta _{p}v+v^{2q-1}\geq
0$ in $D_{\delta }$. Given $M>0$ and $\epsilon >0$, small enough, the set%
\begin{equation*}
\Omega _{\epsilon ,M}=\{x\in \Omega :\epsilon \leq v(x)\leq M\}
\end{equation*}%
is a neighborhood of $\partial \Omega $ \ contained in $D_{\delta }$ (i.e. $%
\Omega _{\epsilon ,M}\subset D$). Then, for any $x_{0}\in \partial \Omega
_{\epsilon ,M},$ we can use a local barrier function $\underline{V}(x)$
based on the expression $c\left\vert x-x_{0}\right\vert ^{\frac{p}{p-2q}}$
over the set $\Omega _{\epsilon ,M}\cap B_{\delta }(x_{0})$, for some $c>0$.
As in the proof of Theorem 2.3 of \cite{Alvarez -Diaz retention}, it is
possible to chose $c>0$ (independent of $h$) such that $\underline{V}(x)$ is
a local subsolution, in the sense that
\begin{equation*}
\left\{
\begin{array}{rl}
-\mu \Delta _{p}\underline{V}+\underline{V}^{2q-1}\leq 0 & \text{in }\Omega
_{\epsilon ,M}\cap B_{\delta }(x_{0}), \\[0.1cm]
\underline{V}\leq v & \text{on }\partial (\Omega _{\epsilon ,M}\cap
B_{\delta }(x_{0})).%
\end{array}%
\right. \text{ }
\end{equation*}%
Thus, by the weak comparison principle $v(x)\geq \underline{V}(x)$ on $%
\Omega _{\epsilon ,M}\cap B_{\delta }(x_{0})$, which implies (\ref{decay q
moderate}) since $\Omega $ is bounded (see an alternative direct proof, for $%
N=1$, in Proposition 1.5 of \cite{Di Interfaces}).

\noindent The proof of i) follows also those type of arguments. Since $q<p$
and $h_{-}(x)\leq \overline{h}_{-}$ on a neighborhood $D_{\delta }$ of $%
\partial \Omega $ we can built a local subsolution $\underline{V}^{\ast }(x)$
on the set $\Omega _{\epsilon ,M}$ (a neighborhood $D_{\delta }$ of $%
\partial \Omega $) such that
\begin{equation*}
-\mu \Delta _{p}\underline{V}^{\ast }+\overline{h}_{-}\underline{V}^{\ast
q-1}\leq 0\text{ in }\Omega _{\epsilon ,M}\cap B_{\delta }(x_{0}),
\end{equation*}%
and the same above arguments apply (leading to the estimate (\ref{decay q
moderate}) since $\Omega $ is bounded) but now building the subsolution by
modifying the function $c\left\vert x-x_{0}\right\vert ^{\frac{p}{p-q}}$
with $c$ depending on $\overline{h}_{-}$.$_{\blacksquare }$

\bigskip

It is useful to study some additional properties satisfied by the
subdifferential $\partial J_{0,q}.$

\begin{lemma}
\label{Lemma compactness}i) $\partial J_{0,q}$ generates a compact semigroup
over $L^{2}(\Omega )$

ii) the resolvent operator $(I+\mu \partial J_{0,q})^{-1}$ leaves invariant
the subspace $L^{\infty }(\Omega )$; i.e. if $h\in L^{\infty }(\Omega )$ and
if $w\in D(\partial J_{0,q}),$ $w\geq 0$, satisfies (\ref{Resolvent equation
2}) then $w\in L^{\infty }(\Omega )$, for any $\mu >0.$
\end{lemma}

\noindent \textit{Proof}. i) Let $\left\{ h_{n}\right\} _{n\in \mathbb{N}}$
be a bounded sequence in $L^{2}(\Omega )$,
\begin{equation*}
\left\Vert h_{n}\right\Vert _{L^{2}(\Omega )}\leq M.
\end{equation*}

\noindent In particular, $h_{n}\rightharpoonup h$ in $L^{2}(\Omega )$ to
some $h\in L^{2}(\Omega ).$ Let $w_{n}\in D(\partial J_{0,q}),$ $w_{n}\geq 0$
be the associated solution of (\ref{Resolvent equation 2}) for any given $%
\mu >0.$ Then, by Lemma \ref{Lemm subdifferential} $v_{n}:=w_{n}^{\frac{1}{q}%
}$ satisfies that $v_{n}\in W_{0}^{1,p}(\Omega )\cap L^{2q}(\Omega )$, $%
\Delta _{p}v_{n},h_{n}(x)v_{n}^{q-1}\in L^{1}(\Omega ),$ $v_{n}$ is positive
and satisfies the sub-homogeneous equation
\begin{equation}
-\mu \Delta _{p}v_{n}+v_{n}^{2q-1}=h_{n}(x)v_{n}^{q-1}\text{ in }\Omega ,
\end{equation}%
in the sense of distributions. By multiplying the equation $w_{n}+\mu
\partial J_{0,q}(w_{n})\ni h_{n}$ by $w_{n},$ from the monotonicity of $%
\partial J_{0,q}$ we get
\begin{equation*}
\left\Vert w_{n}\right\Vert _{L^{2}(\Omega )}\leq M
\end{equation*}%
and so
\begin{equation*}
\left\Vert v_{n}\right\Vert _{L^{2q}(\Omega )}\leq M.
\end{equation*}%
Thus
\begin{equation*}
\left\Vert -\mu \Delta _{p}v_{n}+v_{n}^{2q-1}\right\Vert _{L^{1}(\Omega
)}\leq M^{\prime }
\end{equation*}%
for some $M^{\prime }>0$ (independent on $n$) and thus there exists a
subsequence such that $v_{n}\rightarrow v$ strongly in $L^{1}(\Omega )$ and
weakly in $W^{1,s}(\Omega )$ for any $1\leq s\leq N(p-1)/(N-1)$ (see, e.g.,
\cite{Di} Chapter 4 and its references). By the dominated convergence
Lebesgue theorem $v_{n}^{q}\rightarrow v^{q}$ strongly in $L^{1}(\Omega ).$
Moreover, integrating by parts

\begin{equation*}
\mu \int_{\Omega }|\nabla v|^{p}\mathrm{d}x+\int_{\Omega }v^{2q}\mathrm{d}%
x\leq M^{\prime \prime }
\end{equation*}%
for some $M^{\prime \prime }>0$ and then $v\in W_{0}^{1,p}(\Omega )\cap
L^{2q}(\Omega )$, $\Delta _{p}v,$ $h(x)v^{q-1}\in L^{1}(\Omega )$ (see, e.g.
\cite{Boccardo-Orsina}) and so $w_{n}\rightarrow w$ in $L^{2}(\Omega )$.
Applying the results of \cite{Brezis-commpSemi} (see also Theorem 2.2.2 of
\cite{Vrabie}) we get the conclusion.

\noindent The proof of ii) follows by the Stampacchia iteration method and
it is an obvious modification of Theorem 5.5 of (\cite{Boccardo-Orsina})
(notice that their arguments, for the case $1<q<p$, apply also for this
special purpose to the limit case $q=p$).$_{\blacksquare }$

\bigskip

\begin{remark}
\label{Rem Unique among minima of Euler-Lag}Notice that the functional $%
J_{h,q}$ may have other stationary points different to $w^{1/q}$, with $w$
solution of the resolvent equation (\ref{Resolvent equation 2}). What the
above lemma shows is that the relation $v=w^{1/q}$ gives a uniqueness
criterion for positive solutions of (\ref{Euler-Lagrange}). The positivity
of $v$ is fundamental since it is known that if $\left\vert \left\{ x\in
\Omega :v(x)=0\right\} \right\vert >0$ (which arise, in particular, when $%
h(x)\leq -\overline{h}_{-}<0$ in a neighborhood of $\partial \Omega $ and $%
q<p$ (\cite{Shatzman})) there is multiplicity of nonnegative solutions of (%
\ref{Euler-Lagrange}) (see also \cite{Bandle-Pozio-Tesei}). Nevertheless, if
$q<p$, the uniqueness result applies to `flat solutions' ($i.e.$ positive
solutions such that $u=\frac{\partial u}{\partial n}=0$ on $\Sigma $) (see
\cite{Di-Her-Ilyasov Nonlinear}).When the set $\{x\in \Omega $: $h(x)<0\}$
is \emph{big enough} (or if $\{x\in \Omega $: $h(x)=0\}$ is \emph{big enough}
and $q\in (1,p)$) there are some nonnegative solutions $v$ of (\ref%
{Euler-Lagrange}) which may vanish on some positively measured subset of $%
\Omega $ (and so their support is strictly included in $\overline{\Omega }$)$%
.$ This property (which does not holds when $v=w^{1/q}$ with $w$ solution of
(\ref{Resolvent equation 2})) can be obtained by comparison methods: through
a refined version of \cite{Bensso-Bre-Fri} (see \cite{Di}, \cite{Di. Gaeta}%
), by local energy type methods (\cite{AntDSh libro}), etc.
\end{remark}

\bigskip

\begin{remark}
\label{Rem other powers in resolvent}It is clear that it is possible to
consider \ equations like (\ref{Euler-Lagrange}) with some different
balances between the nonlinear absorption ($v^{2q-1}$) and forcing ($v^{q-1}$%
) terms. Our special case is motivated by the application of the semigroup
theory to the operator $\partial J_{0,q}(w)$ in $L^{2}(\Omega )$.
\end{remark}

\bigskip

\begin{remark}
\label{Remark Lucio} Lemma \ref{Lemm subdifferential} admits many
generalizations dealing with $h\notin L^{2}(\Omega )$ but still with
solutions $v\in W_{0}^{1,p}(\Omega )\cap L^{2q}(\Omega )$\textbf{.} It seems
possible to complement inequality (\ref{L^2 contraction}) by other
inequalities involving different exponents on the norms of the data and the
solutions (see, $e.g.$, \cite{Boccardo-Orsina} and \cite{Palmieri} in the
parabolic framework and Remark \ref{Rem T-accretivity in L^1}).
\end{remark}

\begin{remark}
\label{Rem more general than p-Laplace}It is possible to extend the above
approach by replacing the p-Laplace operator by more general quasilinear
homogeneous operators of the form $\mathrm{div}(a(x,\nabla u))$ with
\begin{equation*}
A(x,t\mathbf{\xi )=}\left\vert t\right\vert ^{p}A(x,\mathbf{\xi )}\text{ for
all }t\in \mathbb{R}\hbox{, } \xi \in \mathbb{R}^{N} \hbox{and a.e. } x\in
\Omega,
\end{equation*}%
where
\begin{equation*}
a(x,\mathbf{\xi })=\frac{1}{p}\partial _{\mathbf{\xi }}A(x,\mathbf{\xi )}
\end{equation*}%
(see \cite{Takac Handbook} and \cite{Girg-Takac}). We point out that the
application of the abstract results of the accretive operators theory allows
also the consideration of this type of diffusion operators (see, $e.g.$,
\cite{Benilan Franco-Japan}).
\end{remark}

\bigskip

\begin{remark}
\label{Rem Strict convexity} A crucial property of the functional $%
J_{0,q}(w) $ is its strict \emph{ray-convexity}: it means that $J_{0,q}(w)$
is strictly convex except for any couple of colinear points$\ w$, $\widehat{w%
}$ with $\widehat{w}$ $=\alpha w$ for some $\alpha \in (0,+\infty )$. That
was used in \cite{Anane}, \cite{Takac-Tello-Ulm} and \cite{Takac Handbook}
to get the uniqueness of nonnegative solutions when $\frac{f_{1}(x,u)}{%
u^{q-1}}$ in (f2) is not strictly decreasing (as it is the case of the first
eigenfunction of the $p$-Laplacian).
\end{remark}

\bigskip

\begin{remark}
\label{Rem p=1 and p=inf} The limit case $p=\infty $ (defined in a suitable
way) can be also considered since, curiously enough, it is an homogeneous
operator of exponent $3$ (see, $e.g.$, \cite{Di Goyo infinite laplacian}).
It is well-known that the other limit case $p=1$ can be also treated as a
subdifferential of a convex function (see $e.g.$, \cite{Andreu-Caseels
D-Mazon}) but the unique choice to apply the reasoning of this paper seems
to be $q=p=1$ and then the results reduce to the well-known case of monotone
perturbations. It would be interesting to know if it is possible to get the
uniqueness of nonnegative solutions of equations\ involving some different
kind of non-monotone sub-homogeneity nonlinear term.
\end{remark}

\section{Selected $L^{s}-$mild solutions, proof of the main theorem and
further remarks}

It is useful to unify the application of abstract results on the associated
Cauchy Problem to the case of the Banach spaces $L^{s}(\Omega )$, for any $%
s\in \lbrack 1,+\infty ].$ For instance, we can define the realizations of
the operator $\partial J_{0,q}$ over the spaces $L^{s}(\Omega )$, for any $%
s\in \lbrack 1,+\infty ]$ as $A_{s}=\overline{\partial J_{0,q}}^{L^{s}}$ in
the sense of graphs over $L^{s}(\Omega )\times L^{s}(\Omega )$: $i.e.$, $%
A_{s}:D(A_{s})\rightarrow \mathcal{P(}L^{s}(\Omega ))$ and $z\in A_{s}(w)$
if and only if there exists $z_{n}\in \partial J_{0,q}(w_{n})$ such that $%
w_{n}\rightarrow w$ and $z_{n}\rightarrow z$ in $L^{s}(\Omega ),$ so that $%
D(A_{s})=$ $\left\{ w\in L^{s}(\Omega ):\exists w_{n}\in L^{2}(\Omega )\text{%
, with }w_{n}^{\frac{1}{q}}\in W_{0}^{1,p}(\Omega )\cap L^{2q}(\Omega )\text{
such that }w_{n}\rightarrow w\text{ in }L^{s}(\Omega )\right\} .$

Then, we consider the Cauchy problem
\begin{equation}
\left\{
\begin{array}{lc}
\dfrac{dw}{dt}+A_{s}w\ni F(t) & \text{in }L^{s}(\Omega ) \\[0.2cm]
w(0)=w_{0}, &
\end{array}%
\right.  \label{Abstract Cauchy}
\end{equation}%
where $w_{0}\in \overline{D(A_{s})}$ and $F\in L^{1}(0,T:L^{s}(\Omega ))$.
In our case, two relevant examples are $A_{2}=\partial J_{0,q}$ and the $%
L^{1}(\Omega )$ operator
\begin{equation*}
\left\{
\begin{array}{l}
AW=-\Delta _{p}W^{m}\text{, for }W\in D(A)\text{, with} \\
D(A)=\{W\in L^{1}(\Omega )\text{, }W^{m}\in W_{0}^{1,1}(\Omega ),\Delta
_{p}W^{m}\in L^{1}(\Omega )\},%
\end{array}%
\right.
\end{equation*}%
given $m>0$ and $p>1.$

We start by recalling the definition of \emph{mild solution} of (\ref%
{Abstract Cauchy}) by particularizing the abstract framework to the case of
the Banach space $X=L^{s}(\Omega )$. The good class of operators to solve (%
\ref{Abstract Cauchy}) is the class of \emph{accretive operators} (resp.
\textit{T-}\emph{accretive operators})\emph{\ }over a Banach space $X:$ i.e.
$A:D(A)\rightarrow \mathcal{P(}X)$ such that
\begin{eqnarray*}
\left\Vert x-\widehat{x}\right\Vert &\leq &\left\Vert x-\widehat{x}+\mu (y-%
\widehat{y})\right\Vert \text{ } \\
\text{(resp. }\left\Vert \left[ x-\widehat{x}\right] _{+}\right\Vert &\leq
&\left\Vert \left[ x-\widehat{x}+\mu (y-\widehat{y})\right] _{+}\right\Vert
\text{)} \\
\text{whenever }\mu &>&0\text{ and }(x,y),(\widehat{x},\widehat{y})\in A.
\end{eqnarray*}%
The operator is called \textit{m-}\emph{accretive }if in addition $R(I+A)=X$%
. For many results and definitions about mild solutions of the Cauchy
Problem for accretive operators in Banach spaces see, e.g., \cite{Barbu},
\cite{Barbu 2010}, \cite{BenCran Pazy}, \cite{Di}, \cite{Veron}, \cite{Evans}
and \cite{Benilan-Wight}. We recall that over any Hilbert space (as $%
L^{2}(\Omega )$) the class of m-T-accetive operators coincides with the
class of maximal T-monotone operators and thus it is possible to apply the
abstract theory presented in Brezis \cite{Brez OMXM}) to problem (\ref%
{Problem parabolic subdifferential J-0q}). The notion of mild solution below
is well defined in both cases: Hilbert and Banach spaces.

\bigskip

\begin{definition}
A function\ $w\in C([0,T]:L^{s}(\Omega ))$ is a $L^{s}-$\emph{mild solution}
of (\ref{Abstract Cauchy}) if for any $\epsilon >0$, there exists a
partition $\{0=t_{0}<t_{1}<...$ $t_{n}\}$ of $[0,t_{n}]$ and there exist two
finite sequences $\{w_{i}\}_{i=0}^{n},$ $\{F_{i}\}_{i=0}^{n}$ in $%
L^{s}(\Omega )$ such that%
\begin{equation*}
\left\{
\begin{array}{l}
\text{(i) }\dfrac{w_{i+1}-w_{i}}{t_{i+1}-t_{i}}+A_{s}w_{i+1}\ni
F_{i+1},\qquad i=0,1,...,n-1 \\[0.2cm]
\text{(ii) }t_{i+1}-t_{i}<\epsilon \\[0.2cm]
\text{(iii) }0\leq T-t_{n}<\epsilon \\[0.2cm]
\text{(iv) }\displaystyle\sum_{i=1}^{n-1}\int_{t_{i}}^{t_{i+1}}\left\Vert
F_{i}-F(t)\right\Vert _{L^{s}(\Omega )}dt<\epsilon ,%
\end{array}%
\right.
\end{equation*}%
and
\begin{equation*}
\left\Vert w_{\epsilon }(t)-w(t)\right\Vert _{L^{s}(\Omega )}\leq \epsilon
\text{ on }[0,t_{n}],
\end{equation*}%
where
\begin{equation*}
w_{\epsilon }(t)=w_{i}\text{ for }t_{i}\leq t<t_{i+1}\text{, }i=0,1,...,n-1.
\end{equation*}
\end{definition}

\begin{definition}
The piecewise constant function $w_{\epsilon }(t)$ defined before is called
an \emph{$\epsilon $-}$L^{s}$-\emph{approximate solution} of (\ref{Abstract
Cauchy}).
\end{definition}

\bigskip

\noindent \textit{Proof of Theorem \ref{Theore main}. }Let us start by
considering the simpler problem $f(x,v)\equiv 0$. Since $w_{0}=u_{0}^{q}\in
D(J_{0,q})\subset \overline{D(J_{0,q})}^{L^{2}}=\overline{D(\partial J_{0,q})%
}^{L^{2}}$, the existence and uniqueness of a mild solution $w\in C([0,%
\widetilde{T}]:L^{2}(\Omega )),$ for any arbitrary $T>0$, is a direct
consequence of the application of the abstract theory (Brezis \cite{Brez
OMXM}) on maximal T-monotone operators in $L^{2}(\Omega )$. Moreover, we
know that $w$ is a weak solution (in the sense of Definition 3.1 of \cite%
{Brez OMXM}): i.e. if we assume $w_{0,n}\in $ $D(\partial J_{0,q})$ and $%
h_{n}\in W^{1,1}(0,\widetilde{T}:L^{2}(\Omega ))$ such that $%
w_{0,n}\rightarrow w_{0}$ in $L^{2}(\Omega )$ and $h_{n}\rightarrow h$ in $%
L^{1}(0,\widetilde{T}:L^{2}(\Omega ))$ then the respective solutions $w_{n}$
satisfy that $w_{n}\rightarrow w$ in $C([0,\widetilde{T}]:L^{2}(\Omega ))$
(see, Theorem 3.4 of \cite{Brez OMXM}). By applying Theorem 3.7 of \cite%
{Brez OMXM} we know that, in fact, $w_{n}$ is a strong solution in the sense
that $w_{n}(t)$ is Lipschitz continuous on $[\delta ,\widetilde{T}]$ for any
$\delta \in (0,\widetilde{T})$ and thus differentiable. Then the associate
problem (\ref{Problem parabolic subdifferential J-0q}) can be written as

\begin{equation*}
\dfrac{dw_{n}}{d\tau }(\tau )-\frac{\Delta _{p}(w_{n}^{{}}(\tau )^{\frac{1}{q%
}})}{w_{n}^{{}}(\tau )^{\frac{q-1}{q}}}=h_{n}(\tau ),
\end{equation*}%
i.e.,
\begin{equation*}
w_{n}^{{}}(\tau )^{\frac{q-1}{q}}\dfrac{dw_{n}}{d\tau }(\tau )-\Delta
_{p}(w_{n}^{{}}(\tau )^{\frac{1}{q}})=h_{n}(\tau )w_{n}^{{}}(\tau )^{\frac{%
q-1}{q}}.
\end{equation*}
If we define $w_{n}^{{}}(\tau )=u_{n}(t)^{q}$ then
\begin{equation*}
w_{n}^{{}}(\tau )^{\frac{q-1}{q}}\dfrac{dw_{n}}{d\tau }(\tau )=\frac{q}{2q-1}%
\dfrac{d(w_{n}^{(2q-1)/q})}{d\tau }(\tau )=\dfrac{d(u_{n}^{2q-1})}{dt}(t)
\end{equation*}%
if
\begin{equation*}
\tau =\frac{q}{2q-1}t.
\end{equation*}%
Obviously we take now $\widetilde{T}=\frac{q}{2q-1}T.$ Notice that $w_{n}\in
C([0,\widetilde{T}]:L^{2}(\Omega ))$ implies that $u_{n}^{q}\in
C([0,T]:L^{2}(\Omega ))$ and thus $u_{n}^{2q-1}\in C([0,T]:L^{2}(\Omega ))$
since $(2q-1)/q>1$ (remember that $q>1$). In addition, for those regular data

\begin{equation*}
\dfrac{d(u_{n}^{2q-1})}{dt}\in \lbrack \widehat{\delta },T]\text{ for any }%
\widehat{\delta }\in (0,T].
\end{equation*}%
Thus, we conclude that $u_{n}(t):=w_{n}(\frac{q}{2q-1}t)^{1/q}$ is a\textit{%
\ }$L^{1}-$\textit{mild positive solution} of $(P_{q})$ on $[0,T]$,
associated to $u_{0,n}:=$ $w_{0,n}{}^{1/q}$ and $h_{n}$ (which the
corresponding unique selected $L^{1}-$\textit{mild positive solution} of $%
(P_{q})$. Finally, as $w_{n}\rightarrow w$ in $C([0,\widetilde{T}%
]:L^{2}(\Omega ))$ we get that $u(t):=w(\frac{q}{2q-1}t)^{1/q}$ is a\textit{%
\ }$L^{1}-$\textit{mild positive solution} of $(P_{q})$ on $[0,T]$,
associated to $u_{0}:=$ $w_{0}{}^{1/q}$ and $h$ since the notion of mild
solution is stable by approximations of the data (see, e.g. Theorem 11. 1 of
\cite{BenCran Pazy}). The rest of conclusions of Theorem \ref{Theore main},
when $f(x,v)\equiv 0$ are a consequence of Lemma \ref{Lemma compactness} and
the T-monotocity of operator $\partial J_{0,q}$ (Lemmas \ref{Prop Benilan}
and \ref{Lemm subdifferential}).

\noindent We consider now the parabolic problem $(P)$ in the general case,
i.e., with a non-homogeneous term $f(x,u)$ satisfying the structural
assumptions (f1)-(f3). We consider now the operator on $L^{2}(\Omega )$
\begin{equation}
Cw=\partial J_{0,q}(w)-\frac{f_{1}(x,w)}{w^{q-1}}  \label{Operator C_1}
\end{equation}%
with $D(C)=D(\partial J_{0,q})$. Since (f1)-(f3) hold and $%
f_{1}(x,w)=f_{1}(w),$ independent of $x$, or $f_{1}(x,w)$ satisfies also
(f4), then the function $E:\Omega \times \lbrack 0,+\infty )\rightarrow
\mathbb{R}$, given by
\begin{equation*}
E(x,w)=-\frac{f_{1}(x,w)}{w^{p-1}}-a_{0}(x)
\end{equation*}%
generates a m--T-accretive operator $L^{2}(\Omega )$ with $E(x,0)=0$. Then,
the operator $C$ is m-T-accretive on $L^{2}(\Omega )$. Moreover, the
Lipschitz function

\begin{equation*}
G(x,w)=-\frac{f_{2}(x,w)}{w^{q-1}}+a_{0}(x)
\end{equation*}%
(of constant $K_{{G}}>0$) generates a Lipschitz operator on $L^{2}(\Omega )$
(of constant $K$ for some $K>0$). Then the operator $C+KI$ is a
m-T-accretive in $L^{2}(\Omega )$ (see, $e.g.$, Chapter 2, Example 2.2 of
\cite{BenCran Pazy}), i.e., $C$ is a $K$-m-T-accretive in $L^{2}(\Omega )$.
So, by the Crandall-Ligget theorem (see, $e.g.$, \cite{Barbu}, and \cite%
{BenCran Pazy}), for any $w_{0}\in \overline{D(\partial J_{0,q})}$ and $h\in
L^{1}(0,T:L^{2}(\Omega ))$ there exists a unique positive $L^{2}-$mild
solution $w\in C([0,T]:L^{2}(\Omega ))$ of the Cauchy Problem
\begin{equation}
\left\{
\begin{array}{lc}
\dfrac{dw}{dt}+\partial J_{0,q}(w)-\dfrac{f_{1}(x,w)}{w^{q-1}}-\dfrac{%
f_{2}(x,w)}{w^{q-1}}\ni h(t) & \text{in }L^{2}(\Omega ) \\[0.2cm]
w(0)=w_{0}, &
\end{array}%
\right.
\end{equation}%
and if $\widehat{w}\in C([0,T]:L^{2}(\Omega ))$ is the $L^{2}-$mild solution
corresponding to the data $\widehat{w}_{0}\in \overline{D(\partial J_{0,q})}$
and $\widehat{h}\in L^{1}(0,T:L^{2}(\Omega ))$ then for any $t\in \lbrack
0,T]$
\begin{equation*}
\begin{array}{ll}
\left\Vert \lbrack w(t)-\widehat{w}(t)]_{+}\right\Vert _{L^{2}(\Omega )} &
\leq e^{Kt}\Vert (w_{0}-\widehat{w}_{0})^{+}\Vert _{L^{2}(\Omega )} \\
& \displaystyle\hspace*{-0.2cm}+\int_{0}^{t}e^{K(t-s)}\left\Vert [h(s)-%
\widehat{h}(s)]_{+}\right\Vert _{L^{2}(\Omega )}ds,%
\end{array}%
\end{equation*}%
(see, $e.g.$, \cite{Barbu 2010} Proposition 4.1 or Theorem 13.1 of \cite%
{BenCran Pazy}). Arguing as before $u(t):=w(\frac{q}{2q-1}t)^{1/q}$ is a%
\textit{\ }$L^{1}-$\textit{mild positive solution} of $(P)$. The proof that $%
u\in L^{\infty }(0,T:L^{\infty }(\Omega ))$ once we assume $h\in L^{\infty
}(0,T:L^{\infty }(\Omega ))$ and $u_{0}\in L^{\infty }(\Omega )$ is a
consequence of Lemma \ref{Lemma compactness} (which implies the compactness
of the semigroup generated by operator $\partial J_{0,q}(w)-\dfrac{f_{1}(x,w)%
}{w^{q-1}}-\dfrac{f_{2}(x,w)}{w^{q-1}}$) and the abstract invariant results
presented in Theorem 2.4.1 of Vrabie \cite{Vrabie} (see also \cite{Di Vrabie}%
), which ends the proof of Theorem \ref{Theore main} $_{\blacksquare }$

\bigskip

\begin{remark}
\label{Rem Hipo dato inicial}In fact, the existence and uniqueness of a $%
L^{2}-$mild positive solution of problem (\ref{Problem parabolic
subdifferential J-0q}) can be assured in the more general case of $w_{0}\in
\overline{{\mathcal{D}}(\partial J_{0,q})}$. Notice that if $w_{0}\in
\overline{{\mathcal{D}}(\partial J_{0,q})}$ the selected $L^{1}-$mild
positive solution $u$ of $(P_{q})$ such that $w(\frac{q}{2q-1}t)=u(t)^{q},$
with $w(t)$ the corresponding $L^{2}-$mild positive solution of problem (\ref%
{Problem parabolic subdifferential J-0q}) satisfies (in some sense) the
decay estimates given in Lemma \ref{Lemm subdifferential} since they are
obtained through the implicit Euler scheme given in the definition of mild
solution. As a matter of fact, if $w(t_{0})\in {\mathcal{D}}(\partial
J_{0,q})$ for some $t_{0}\in \lbrack 0,T]$, i.e. $\partial
J_{0,q}(w(t_{0}))\ni h(t_{0})$ for some $h(t_{0})\in L^{2}(\Omega )$ then $%
-\Delta
_{p}v(t_{0})+h(t_{0})_{-}(x)v(t_{0})^{q-1}=h(t_{0})_{+}(x)v(t_{0})^{q-1}$
and necessarily we get the estimates iii) and iv) of \ref{Lemm
subdifferential} for $v(t_{0})$. We also point out that some uniqueness
results for suitable sublinear parabolic problems, when $u_{0}(x)\geq
Cd(x,\partial \Omega ),$ can be found in \cite{Cazenave-Dik-Escobedo}, \cite%
{Giacomoni-Sauvy-Shmarev}, \cite{Dao-Diaz}, \cite{Di-Her-Ilyasov}, \cite%
{Diaz-Giacomoni} (see also their references to previous works in this
direction). Curiously enough such type of assumptions also lead to the
uniqueness of solutions in the case of equations with multivalued right hand
side terms as problem $(P_{H})$ (see \cite{Feireisel}, \cite{Di-Tello})
which until now required completely different ideas.
\end{remark}

\bigskip

\begin{remark}
We point out that selected $L^{1}-$mild positive solution $u$ satisfies some
extra regularity properties due to the subdifferential of $J_{0,q}$ involved
in the equation. See also some variational type techniques applied to the
case $p=2$ in \cite{Palmieri} and the general approach (also for $p=2$)
presented to some related problems in \cite{Brezis-Cacenave}, \cite{Brez
Cazenave libro}.
\end{remark}

\bigskip

\begin{remark}
It seems possible to make a sharper study of the regularity of the solution
of the equation $-\mu \Delta _{p}v+v^{2q-1}=h(x)v^{q-1}$, but we shall not
enter into the maximum of its generality here. For instance, when $p=2$ such
equation becomes a Schr\"{o}dinger equation with a potential $h(x)$ (and a
nonlinear perturbation term $v^{3}$) and so it is possible to consider
potentials $h(x)$ with a singular behavior near $\partial \Omega $ (and in
other subregions of $\Omega $) which goes beyond $L^{1}(\Omega )$ (see, $%
e.g. $, \cite{Ben-Petra-absorptions}, \cite{Ponce}, \cite{Di Go Rako},\cite%
{Orsina ponce}, \cite{Di Go Vazq} and its many references). For the special
case of $q=p\neq 2$ singular potentials were considered in \cite%
{Marcus-Shafir}, \cite{Pincho-Tertikas-Titore}, \cite{DiHernMance2} and in
many other papers.
\end{remark}

\bigskip

\begin{remark}
\label{Rem Caso R^N} The main result of this paper may be also proved when
we replace the open bounded set $\Omega $ by the whole space $\mathbb{R}%
^{N}. $ The \emph{Diaz-Sa\'{a} inequality} (and the generalized Picone
inequality) was obtained in \cite{Chaib} (respectively in \cite{Dat}). We do
no want to enter into details here but the arguments of truncating the
domain, generate the associate problems on an expansive sequence of domains $%
\Omega _{n}$ and then to get the solution as limit of the solutions of the
corresponding problems on $\Omega _{n}$ can be applied as in Brezis and
Kamin \cite{Br-Kamin} (see also \cite{Di Oleinik}). The assumptions made on
functions $f_{i}$ allow to get some similar estimates to (\ref{Estimate
dependence parabolic}) to solutions of several quasilinear formulations
(see, \cite{Pablo-Vazquez}, and \cite{Hui}) and, in particular, to solutions
of the associated to the KPP equation as in the papers \cite{Br-Kamin}, \cite%
{Di-Kamin}, \cite{Audrito-Vazquez} and \cite{Audruito.Vazquez slow}).
\end{remark}

\bigskip

\begin{remark}
\label{rem singular equ}As mentioned before, the assumptions on $f_{1}(x,u)$
allow the consideration of some singular terms: see, $e.g.$, \cite%
{Badra-Giacomoni}, \cite{Brou-Giacomoni}, \cite{Ghergu-Radulescu}, \cite%
{Diaz-Giacomoni} and the surveys \cite{Hernandez-Mancebo} and \cite%
{Ghergu-Radulescu libro singular}. The assumption of the type $\frac{%
f_{2}(x,u)}{u^{q-1}}$ globally Lipschitz continuous in $u\in (0,+\infty )$
was used for other purposes in previous works in the literature (see, $e.g.$%
, \cite{Cuesta -Takak}).
\end{remark}

\bigskip

\begin{remark}
\label{rEM nEUMANN TYE boundary cnditions} It seems possible to get similar
results to positive solutions of Neumann type boundary conditions once that
the homogeneity of the boundary condition is compatible with the one of the
doubly nonlinear problem $(P)$ (see, $e.g.$, \cite{Benilan Franco-Japan},
\cite{Alt-Luckhaus}, \cite{Bandle-Pozio-Tesei} and \cite{Andreucci} among
many other possible references).
\end{remark}

\bigskip

We point out that, obviously, the function $u_{\infty }(x)\equiv 0$ in $%
\Omega $ is a trivial solution of the stationary problem. Here we are
interested on nonnegative solutions of problem $(P)$ (and its implicit time
discretization). The following result shows that the asymptotic behavior, as
$t\rightarrow +\infty $, is very different according $q\in (1,\frac{p}{2})$
and $p>2$ than in the case $q\in (\frac{p}{2},p].$

We will prove that, in fact, if $q\in (1,\frac{p}{2})$ and $p>2,f(x,u)\equiv
0,~h\in L^{1}(0,T:L^{2}(\Omega )),$ $h\geq 0$ and $u_{0}\gvertneqq 0$ then
there is no extinction in finite time, so that $\Vert u^{2q-1}(t)\Vert
_{L^{2}(\Omega )}>0$ for any $t>0.$ The situation is different if $q\in (%
\frac{p}{2},p]$ since, at least for $f(x,u)\equiv 0$ and $h\leq 0$, there is
a finite extinction time $T_{e}>0$, such that $w(t)\equiv 0,$ in $\Omega $,
for any $t\geq T_{e}$. In that case, we understand that the selected
solution $v(t)$ of $(P)$ also extinguishes in $\Omega $ after $T_{e}.$

\bigskip

\begin{theorem}
\label{Theorem extinction copy(1)} a) Assume $q\in (1,\frac{p}{2})$ and $%
p>2,f(x,u)\equiv 0,~h\in L_{loc}^{1}(0,+\infty :L^{2}(\Omega )),$ $h\geq 0$
and $u_{0}\gvertneqq 0$ satisfying (\ref{Hypo u_0}). Then the selected $%
L^{1}-$mild positive solution $u$ of $(P)$ satisfies that%
\begin{equation*}
\Vert u^{q}(t)\Vert _{L^{2}(\Omega )}\geq \frac{1}{%
(c_{1}t+c_{2})^{(q-1)/(p+q-2)}}
\end{equation*}
for any $t>0$, for some positive constants $c_{1}$ and $c_{2}.$

b) Assume $q\in (\frac{p}{2},p]$, $f(x,v)\equiv 0$ and $h\in
L_{loc}^{1}(0,+\infty :L^{2}(\Omega ))$ such that $h\leq 0.$Then there is a
finite extinction time $T_{e}>0$, such that the selected solution $u(t)$ of $%
(P)$ extinguishes in $\Omega $ after $T_{e},$ i.e., $u(t)=u_{\infty
}(x)\equiv 0,$ in $\Omega $, for any $t\geq T_{e}$.
\end{theorem}

\noindent \textit{Proof}. Since $h\geq 0$, from the comparison estimate (\ref%
{Estimate dependence parabolic}) we deduce that $u\geq \underline{U}$ with $%
\underline{U}$ the unique solution of the problem%
\begin{equation*}
(P_{0})\qquad \left\{
\begin{array}{rclr}
\partial _{t}(\underline{U}^{2q-1})-\Delta _{p}\underline{U} & = & 0 & \text{
in }Q_{T}, \\[0.1cm]
\underline{U} & = & 0\  & \text{ on }\Sigma , \\[0.1cm]
\underline{U}^{q}(0,.) & = & u_{0}^{q}(.) & \text{ on }\Omega .%
\end{array}%
\right.
\end{equation*}%
Moreover, as indicated in Theorem \ref{Theore main}, we know that if $%
\underline{U}(t):=\underline{W}(\frac{q}{2q-1}t)^{1/q}$ then $\underline{W}$
satisfies of the problem
\begin{equation}
\left\{
\begin{array}{lc}
\dfrac{d\underline{W}}{dt}+\partial J_{0,q}(\underline{W})\ni 0 & \text{in }%
L^{2}(\Omega ) \\[0.2cm]
\underline{W}(0)=u_{0}. &
\end{array}%
\right.
\end{equation}%
In addition, the operator $\partial J_{0,q}(\underline{W})$ is formally
given by $\frac{\Delta _{p}(w^{\frac{1}{q}})}{w^{\frac{q-1}{q}}}$ and thus
it is homogeneous of exponent $\theta =(p-q)/q$, in the sense that
\begin{equation*}
\partial J_{0,q}(r\underline{W})=r^{\theta }\partial J_{0,q}(r\underline{W})%
\text{ for any }r\geq 0\text{ and }\underline{W}\in D(\partial J_{0,q}).
\end{equation*}%
Then, since $q\in (1,\frac{p}{2})$ and $p>2$ implies that $\theta >1$,
applying Theorem 1.1 of \cite{Alikakos-Rostamian} we get that
\begin{equation*}
\Vert \underline{U}^{q}(t)\Vert _{L^{2}(\Omega )}\geq \frac{1}{%
(c_{1}t+c_{2})^{(q-1)/(p+q-2)}}\text{ for any }t>0,
\end{equation*}%
for some positive constants $c_{1}$ and $c_{2}$, and then the conclusion
holds since $U\geq \underline{U}.$

\noindent b) We consider, again the solution $\underline{U}$ of $(P_{0}).$
Now $0\leq u\leq \underline{U}$ \ and since, in this case, the homogeneity
exponent of $\partial J_{0,q}(\underline{W})$ is $\theta <1$ the conclusion
results of the application of Corollary 1 of \cite{Belaud-Diaz}. $%
_{\blacksquare }$

\bigskip

\begin{remark}
\label{Rem sytems and higher order} \ Systems involving sub-homogeneous
terms have been extensively considered in the literature: see, $e.g.$, \cite%
{Fle-Hern Thelin systems}, \cite{Flekinger Gossez Hernandez}, \cite{Chaib}
and its references. It would be interesting to apply the assumptions of the
general framework in this paper to the case of systems. In the case of \emph{%
higher order equations with sub-homogeneous} terms the T-accretivity in $%
L^{p}$ fails but I conjecture that the $L^{2}-$contraction continuous
dependence still holds for certain homogeneous higher order operators (as
for instance those considered in \cite{Bernis} and \cite{Alt}).
\end{remark}

\bigskip

\begin{remark}
\label{Rem T-accretivity in L^1}As mentioned before (see Remark \ref{Remark
correccion k>0}) a stronger property on the convex funcional $J$ may lead to
the accretivity in $L^{1}$ and in $L^{\infty }$ of the realization over
these spaces of the subdifferential operator $\partial J$. Although we are
not able to check the stronger property (\ref{tau-contraction for any k}) in
the special case of functional $J_{0,q}$ it is possible to get some
continuity dependence inequalities for solutions of the equation $w+\mu
\partial J_{0,q}(w)\ni h$, for any $\mu >0$, which keep some resemblances
with the inequalities expressing the $L^{1}$ and $L^{\infty }$ T-accretivity
for the realization of the operator $\partial J_{0,q}(w)$ over those spaces
(some related techniques can be found in Brezis and Kamin \cite{Br-Kamin}
and \cite{Charro-Peral}).
\end{remark}

\bigskip

\begin{acknowledgement}
It is a great pleasure to thank the many discussions with Jacques Giacomoni
on a very preliminary version of this paper. In particular, he showed me how
to prove that the operator $\partial J_{0,q}$ is m-T-accretive in $%
L^{2}(\Omega )$ by using the Picone inequality instead the convexity of $%
J_{0,q}$. I also thank Lucio Boccardo for several comments (in particular on
Remark \ref{Remark Lucio}) and Gregorio D\'{\i}az, David G\'{o}mez-Castro,
Jes\'{u}s Hern\'{a}ndez, Jean Michel Rakotoson and Laurent Veron for some
useful conversations. The research was partially supported by the project
ref. MTM2017-85449-P of the DGISPI (Spain) and the Research Group MOMAT
(Ref. 910480) of the UCM.
\end{acknowledgement}

\end{document}